\documentclass[11pt,reqno]{article}

\usepackage{amsfonts}
\usepackage{amsmath}
\usepackage{amsthm}
\usepackage{amssymb}
\usepackage{mathrsfs}
\usepackage{amstext}

\usepackage[dvips]{graphicx}
\usepackage{subfigure}
\usepackage{algorithm}
\usepackage{algorithmic}

\usepackage{epstopdf}

\usepackage{url}

\font\msbm=msbm10

\numberwithin{equation}{section}

\theoremstyle{plain}
\newtheorem{Theorem}{Theorem}[section]
\newtheorem{lemma}[Theorem]{Lemma}
\newtheorem{corollary}[Theorem]{Corollary}

\newtheorem{definition}{Definition}[section]

\theoremstyle{remark}
\newtheorem{remark}{Remark}[section]

\def\mathbb#1{\hbox{\msbm{#1}}}
\newcommand{\N}{{\mathbb{N}}}

\newcommand{\Z}{{\mathbb{Z}}}
\newcommand{\C}{{\mathbb{C}}}

\renewcommand{\P}{{\mathbb{P}}}
\newcommand{\E}{{\mathbb{E}}}

\newcommand{\on}{{|\!|\!|}}

\newcommand{\beq}{\begin{eqnarray}}
\newcommand{\eeq}{\end{eqnarray}}
\newcommand{\beqn}{\begin{eqnarray*}}
\newcommand{\eeqn}{\end{eqnarray*}}

\newcommand{\ff}{{\mathrm f}}
\newcommand{\F}{{\cal F}}

\newcommand{\supp}{\operatorname{supp}}

\newcommand{\spn}{\operatorname{span}}

\renewcommand{\qed}{\rule{2.5mm}{2.5mm}}

\begin{document}

\title{Stability Results for Random Sampling of Sparse Trigonometric Polynomials}
\author{Holger Rauhut}
\date{September 18, 2006; revised April 4, 2008}
%NuHAG, Faculty of Mathematics, University of Vienna\\
%Nordbergstrasse 15, A-1090 Wien, Austria\\
%rauhut@ma.tum.de

\maketitle

\begin{abstract}
Recently, it has been observed that a sparse trigonometric polynomial, 
i.e.~having only a small number of non-zero coefficients, can be
reconstructed exactly from a small number of random samples using
Basis Pursuit (BP) or Orthogonal Matching Pursuit (OMP). 
In the present article it is shown that recovery by
a BP variant is stable
under perturbation of the samples values by noise. A 
similar partial result for OMP is provided.
For BP in addition, the stability result
is extended to (non-sparse) trigonometric polynomials 
that can be well-approximated by sparse
ones. The theoretical findings are illustrated by numerical experiments.
\end{abstract}
\vspace{0.5cm}


%94A20: Sampling Theory
%42A05: Trigonometric Polynomials
%15A52: Random matrices
%05A18: Enumerative Combinatorics: Partitions of Sets
%90C05: linear programming
%90C25: convex programming

\section{Introduction}

Over the recent years compressed sensing 
has become a rapidly developing research field, 
see e.g.~\cite{babadusawa05,carota06,codade06,do06-2,ru06-1,dagistzo05}. 
In their seminal papers
\cite{carota06,carota06-1,cata06} 
Cand{\`e}s, Romberg and Tao observed that it is 
possible to recover sparse vectors, 
i.e., having only few non-vanishing coefficients, from a 
number of measurements that is small compared to the ambient dimension
of the vector. As reconstruction method they promoted $\ell_1$-minimization, 
also refered to as Basis Pursuit (BP) \cite{chdosa99}. 
Their results apply in particular to recovery
of a sparse vector from (random) samples of its discrete Fourier transform.
In \cite{ra05-7} the author extended their result to the situation
where samples of the corresponding trigonometric polynomial are taken 
at random from the uniform (continuous) distribution on the cube, i.e.,
the samples are chosen ``off the grid''. 

Another line of research suggests Orthogonal Matching Pursuit (OMP)
as recovery method \cite{gitr05,kura06,tr04}. 
This is a greedy algorithm which is significantly faster than BP in practice. 
Partial results in \cite{kura06} indicate
that also OMP is able to recover
a sparse trigonometric polynomial from few random samples. 
Moreover, numerical
experiments suggest that OMP usually
has a slightly higher probability of recovery 
success than BP - although BP has some theoretical advantages.

In practice, it is important that recovery methods are stable in the presence
of noise on the measurements. Cand{\`e}s et al.~showed 
in \cite{carota06-1} that
(a variant of) BP is indeed stable under a certain
condition on the measurement matrix involving the so called
restricted isometry constants. An estimation of these
constants for the measurement matrix corresponding
to random samples of the discrete Fourier transform
was provided in \cite{cata06} and \cite{ru06-1}. 
In the present article we extend this estimate
to the case of random samples at uniformly distributed points on
the cube $[0,2\pi]^d$.

We further provide partial results indicating that also OMP is stable 
under perturbation of the measurements by noise. 
Finally, numerical
experiments reveal that the average reconstruction error of OMP is
usually smaller than for (the variant of) BP in the presence of noise.

After the first submission of this manuscript, variants of OMP -- 
Regularized Orthogonal Matching Pursuit (ROMP) \cite{neve07,neve07-1}
and CoSaMP \cite{netr08} -- were introduced, that achieve similar theoretical
recovery and stability guarantees as Basis Pursuit and are even slightly 
faster than OMP. Since the analysis of these algorithms 
is based on the restricted isometry constants
our estimates for the Fourier type measurement matrix are 
useful for the analysis of ROMP and CoSaMP as well.

The paper is organized as follows. Section 2 gives some background on prior
work, introduces notation and describes our problem. In Section 3 we present 
our main results concerning stability of a variant of BP, 
while Section 4 states stability theorems
for OMP. Section 5 presents the proofs for BP, and Section 6 deals
with the ones for OMP. The numerical experiments are detailed in Section 7.
Finally, we conclude in Section 8 with a discussion.

\section{Prior Work and Problem Statement}
\label{Sec_Prior}

For some finite subset $\Gamma \subset \Z^d$, $d \in \N$, we let $\Pi_\Gamma$
denote the space of all trigonometric polynomials in dimension $d$ whose
coefficients are supported on $\Gamma$.
An element $f$ of $\Pi_\Gamma$ is of the form 
$f(x) = \sum_{k \in \Gamma} c_k e^{i k\cdot x}$, $x \in [0,2\pi]^d$, 
with Fourier coefficients $c_k \in \C$.
The dimension of $\Pi_\Gamma$ will be denoted by $D:= |\Gamma|$.
One may imagine $\Gamma = \{-q,-q+1,\hdots,q-1,q\}^d$, but actually arbitrary
sets $\Gamma$ are possible.

We will mainly deal with ``sparse'' trigonometric polynomials, i.e., we assume
that the sequence of coefficients $c_k$ is supported 
only on a small set $T \subset \Gamma$.
However, a priori nothing is known about $T$ apart from a maximum size.
Thus, it is useful to introduce the (nonlinear) 
set $\Pi_\Gamma(M) \subset \Pi_\Gamma$ of
all trigonometric polynomials whose Fourier coefficients are supported on a
set  $T \subset \Gamma$ satisfying $|T| \leq M$, 
$\Pi_\Gamma(M)\,=\, \bigcup_{T\subset \Gamma, |T|\leq M} \Pi_T$.

Our aim is to reconstruct an element $f \in \Pi_\Gamma(M)$ from
sample values $f(x_1),\hdots, f(x_N)$, where the number $N$ of sampling points
$x_1, \hdots, x_N \in [0,2\pi]^d$
is small compared to the dimension $D$ (but, of course, larger than the 
sparsity $M$). As suggested by \cite{carota06,cata06,gitr05,kura06,ra05-7}
we will study the behaviour of two
reconstruction methods: Basis Pursuit (BP) 
and Orthogonal Matching Pursuit (OMP).

BP was much promoted by Donoho and his coworkers, 
see e.g.~\cite{chdosa99,dota05}. 
It consists in solving the following $\ell^1$-minimization problem
\begin{equation}\label{BP}
\min \|(d_k)\|_1 = \sum_{k \in \Gamma} |d_k| \quad \mbox{ subject to } 
\quad \sum_{k \in \Gamma } 
d_k e^{2\pi i k\cdot x_j} = f\left(x_j\right),\quad j=1,\hdots,N.
\end{equation}
This task can be performed with convex optimization techniques \cite{bova04}.
Recently, much effort has been dedicated to the development of fast
algorithms specialized to $\ell_1$-minimization, see e.g.\ \cite{dafoloXX,finowr07,bogokikslu07}.

OMP is a greedy algorithm 
\cite{mazh93,tr04}, which selects
a new element of the support set $T$ in each step, 
see Algorithm \ref{algo:omp}. 
Its precise formulation uses the following notation. 
Let $X=(x_1,\hdots,x_N)$ be the sequence of sampling points.  
We denote by $\F_X$ the $N \times D$ matrix with
entries
\begin{equation}\label{def_FX}
  (\F_X)_{j,k} \,=\, e^{ik \cdot x_j}, \quad 1\leq j \leq N,\, k \in \Gamma.
\end{equation}
Then clearly, $f(x_j) = (\F_X c)_j$ if $c$ is the vector of Fourier
coefficients of $f$. Let $\phi_k$ denote the $k$-th column of $\F_X$, i.e.,
$\phi_k \,=\, \left(e^{ik\cdot x_\ell}\right)_{\ell=1}^N$.
The restriction
of $\F_X$ to the columns indexed by $T$ is denoted by $\F_{TX}$.
Furthermore, let $\langle \cdot, \cdot \rangle$ denote the usual Euclidean
scalar product and $\|\cdot\|_2$ the associated norm.
We have $\| \phi_k \|_2 = \sqrt{N}$ for all $k \in \Gamma$, i.e., all the
columns of $\F_X$ have the same $\ell^2$-norm. For details on the 
implementation of OMP we refer to \cite{kura06}. We only note that the fast
Fourier transform (FFT) or the non-equispaced fast Fourier transform (NFFT),
see e.g.~\cite{postta01} and the references therein, 
can be used for speed-ups of OMP.

\begin{algorithm}[t]
  \caption{OMP\label{algo:omp}}
  \begin{tabular}{ll}
    Input:    & sampling set $X\subset [0,2\pi]^d$, sampling vector
                $\ff:=(f(x_j))_{j=1}^N$, set $\Gamma\subset\Z^d$.\\
    Optional: & maximum allowed sparsity $M$ and/or residual tolerance
                $\varepsilon$.\\[0.5ex]
  \end{tabular}
  \begin{algorithmic}[1]
    \STATE Set $s=0$, the residual vector $r_0 = \ff$, and the index set
           $T_0 = \emptyset$.
    \REPEAT
    \STATE Set $s=s+1$.
    \STATE Find $k_s=\arg\max_{k\in\Gamma} |\langle r_{s-1},\phi_k\rangle|$
    and augment $T_s = T_{s-1} \cup \{k_s\}$.
    \STATE Project onto $\spn\{ \phi_k, k\in T_s\}$ by solving the least squares problem
    \begin{equation*}
      \left\|\F_{T_s X} d_s - \ff\right\|_2\stackrel{d_s}{\rightarrow}\min.
    \end{equation*}
    \STATE Compute the new residual $r_s =  \ff - \F_{T_s X} d_s$.
    \UNTIL{$s=M$ or $\|r_s\|\le\varepsilon$}
    \STATE Set $T = T_s$, the non-zeros of the vector $c$ are given by
    $(c_k)_{k\in T}=d_s$.
  \end{algorithmic}

  \begin{tabular}{ll}
    Output:   & vector of coefficients $(c_k)_{k\in \Gamma}$ and its support
                $T$.
  \end{tabular}
\end{algorithm}

Since it seems to be very hard to come up with deterministic recovery results we model
the sampling points $x_1,\hdots,x_N$ as random variables. To this end we use two
probability models. 
\begin{itemize}
\item[(1)] The sampling points $x_1,\hdots,x_N$ are independent 
random variables having the uniform distribution on the cube $[0,2\pi]^d$.
\item[(2)] The sampling points $x_1,\hdots,x_N$ are independent
random variables having the uniform distribution on the grid 
$\{0,\frac{2\pi}{q},\hdots,2\pi\frac{q-1}{q}\}^d$. Here, it is implicitly
assumed that $\Gamma \subseteq \{0,1,\hdots,q-1\}^d$.
\end{itemize}
Model (1) will also be refered to as the continuous model, while the second
will be called ``discrete''. Observe that with model (2) it might happen with 
non-zero probability that some sampling points are selected more than once.
To overcome
this problem one might also choose the sampling set uniformly
at random among all subsets of the grid $\{0,\frac{2\pi}{q},\hdots,2\pi\frac{q-1}{q}\}^d$
of size $N$. This model was actually used in \cite{carota06,cata06,ru06-1}.
However, for technical reasons we work with the model (2) here. 
Intuitively, moving from model (2) to its variant should actually improve the situation
since always a maximum of information is used.

In \cite{ra05-7} it was proven that BP 
is able to recover
a sparse trigonometric polynomial from a rather small number of
sample values.

\begin{Theorem}\label{thm_old} Let $T \subset \Gamma$ with $|T| \leq M$. Choose $x_1,\hdots,x_N$
be random variables according to the probability models (1) or (2). Assume that
\begin{equation}\label{N:cond}
N \geq C M \log(D/\epsilon).
\end{equation}
Then with probability at least $1-\epsilon$ both BP and OMP recover exactly
all $f \in \Pi_\Gamma(M)$ with coefficients supported on $T$ from
the sample values $f(x_j), j=1,\hdots,N$. 
The constant $C$ is absolute.
\end{Theorem}

The above theorem is non-uniform in the sense that for a single sampling set
$X$ recovery is guaranteed only for the given support set $T$ (but for all
Fourier coefficients supported on $T$). By Theorem  
\ref{thm_BP} to be shown later it follows that this drawback can 
be removed, i.e., recovery by BP can be made fully uniform
by introducing additional $\log$ factors to condition (\ref{N:cond}).

Recovery by OMP was studied theoretically and numerically 
in \cite{kura06}, although the theoretical results are only partial so far.
At least the first step of OMP could be analyzed:

\begin{Theorem}\label{thm:OMP_old} Let $f\in \Pi_\Gamma(M)$ with coefficients
supported on $T$. Choose random sampling points $x_1,\hdots,x_N$ according
to one of our two probability models. If
\[
N \geq C M \log(D/\epsilon)
\]
then with probability at least $1-\epsilon$ OMP 
selects an element
of the true support $T$ in the first iteration.
\end{Theorem}
The numerical experiments conducted in \cite{kura06} suggest that also
the further steps of OMP select elements of the true support $T$, so that
after $M$ steps the correct polynomial $f$ is recovered. However, starting with
the second step the theoretical analysis seems to be quite difficult due
to subtle stochastic dependency issues. 

We note that the above theorem is non-uniform in the sense that the 
success probability is valid for the given polynomial, but it does not
state that with high probability a single sampling set $\{x_1,\hdots,x_N\}$
is good for {\em all} sparse trigonometric polynomials. Such a uniform result
was also provided in \cite{kura06}, which actually analyzes the full application
of OMP, but requires significantly more samples.

\begin{Theorem}\label{thm:OMP_uniform}  
Let $X=(x_1,\hdots,x_N)$ be chosen according
to the continuous probability model (1) or the discrete model (2).
Suppose that
\begin{equation}\label{cond:coherenceN}
N \geq C (2M-1)^2 \log(4D'/\epsilon),
\end{equation}
where $D' := \#\{j-k: j,k \in \Gamma, j\neq k\} \leq D^2$.
Then with probability at least $1-\epsilon$ 
OMP recovers {\it every} $f \in \Pi_\Gamma(M)$. The constant satisfies 
$C \leq 4 +\frac{4}{3\sqrt{2}} \approx 4.94$. In case of 
the continuous probability model it can be improved to $C = 4/3$.
\end{Theorem}

The above result is based on analysis of the coherence, see also below.
It seems that condition (\ref{cond:coherenceN}) is actually
optimal up to perhaps the constant $C$ and the $\log$-factor if one
requires uniformity, i.e., recovery of {\em all} sparse trigonometric
polynomials in $\Pi_\Gamma(M)$ from a single sampling set $X$, see \cite{ra07}.
In this regard, BP and OMP seem to be crucially different.
BP can give a uniform guarantee if the number of samples $N$
scales linearly in the sparsity $M$ (ignoring $\log$-factors), see 
Theorem \ref{thm_old}, Theorem \ref{thm_BP} below and 
e.g.~\cite{cata06,ru06-1}, while OMP can give at most a non-uniform guarantee
in this range, compare also \cite[Section 7]{do06} and \cite{gitr05}.

\bigskip

In this article we treat 
the question whether recovery by BP and OMP is stable if the sample
values $f(x_j)$ are perturbed by noise. Additionally, for BP
we consider also the case that $f$ is not sparse in a strict sense, but can be
well approximated by a sparse trigonometric polynomial.

\begin{figure}%[h]
  \centering
  \subfigure[Trig. polynomial and (noisy) samples.]   
  {\includegraphics[width=0.45\textwidth]{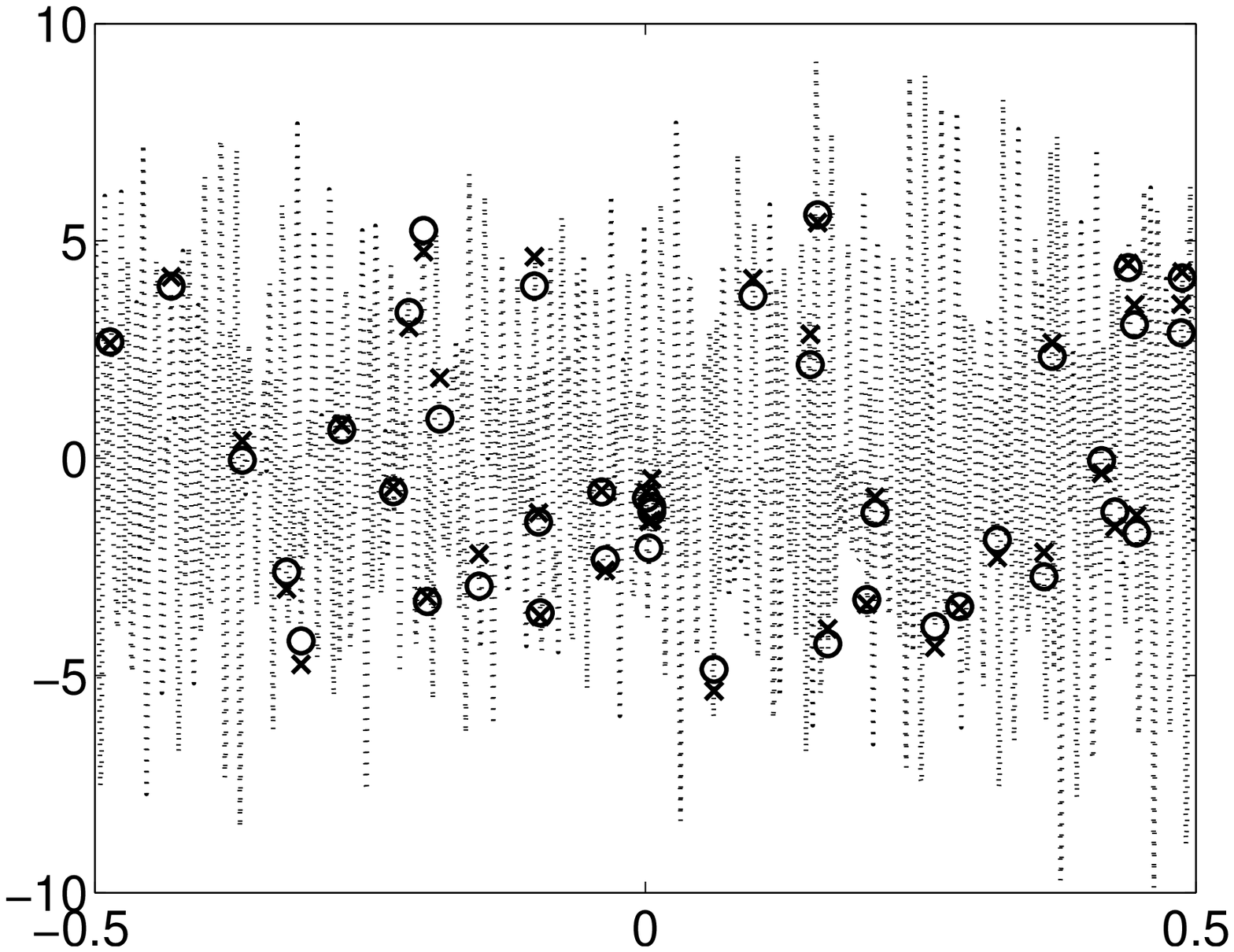}}\hfill
  \subfigure[True and recovered coefficients.]
  {\includegraphics[width=0.45\textwidth]{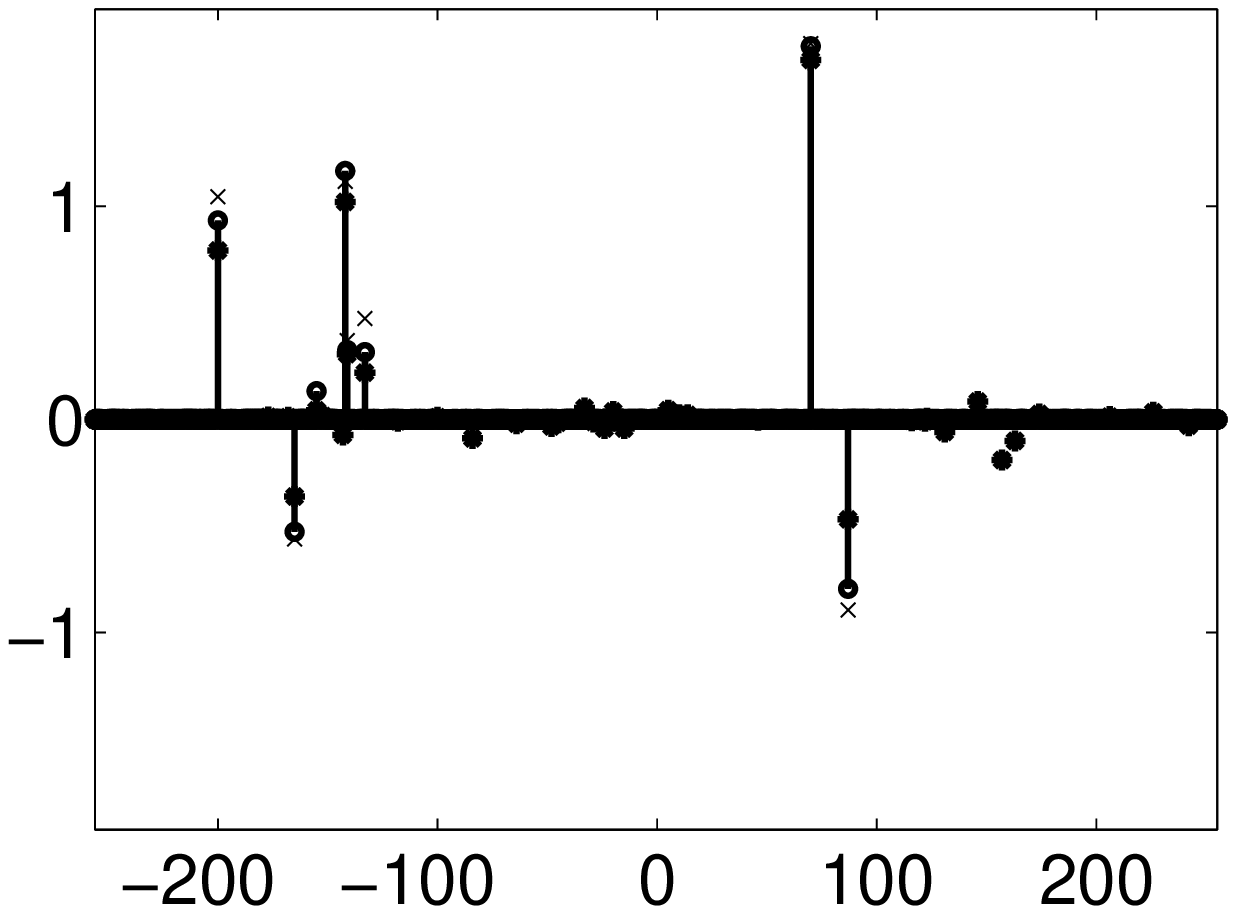}}\\
  \caption{Left: Trigonometric polynomial (real part) of sparsity $M=8$ and
    $N=40$ samples ($\circ$).
    The samples are disturbed by noise $\eta$ 
    with $\|\eta\|_2 = 4$ ($\times$). 
    Right: True coefficients ($\circ$), reconstruction by BP variant (\ref{P2}) ($*$),
    reconstruction by OMP ($\times$).\label{fig:1}}
\end{figure}

In mathematical terms we assume that we observe the vector
\[ 
y \,=\, (f(x_j))_{j=1}^N + \eta \,=\, \F_X c + \eta,
\] 
rather than $(f(x_j)) = \F_X c$,
where the noise $\eta$ satisfies 
$\|\eta\|_2 = (\sum_{\ell=1}^N \eta_\ell^2)^{1/2} \leq \sigma$ for some $\sigma \geq 0$.
We will investigate whether the difference between the original
coefficient vector and the one reconstructed by OMP or BP is small. 
For OMP we additionally ask whether the correct support set is recovered.
Figure \ref{fig:1} provides a first illustration by showing 
an example of a reconstruction by the BP variant
(\ref{P2}) and OMP from noisy samples.

In the sequel, $\|\cdot\|_{p\to q}$ will denote the operator norm from 
the sequence space
$\ell^p$ into $\ell^q$ (on some index set), $\lfloor x \rfloor$ is the largest
integer smaller or equal to $x$. Furthermore, $C$ will always 
denote a generic constant, whose value might be
different in each occurence.

\section{Basis Pursuit}

In the presence of noise it is useful to consider
a slight variant of Basis Pursuit. Indeed, in \cite{carota06-1} it is suggested
to minimize the $\ell_1$-norm
of the coefficient vector $c$ subject to the constraint that the residual 
error satisfies $\|\F_X c - y\|_2 \leq \sigma$, i.e., we solve
\begin{equation}\label{P2}
\min \|c\|_1 \quad \mbox{ subject to } \|\F_X c - y \|_2 \leq \sigma.
\end{equation}
Again this problem can be solved by 
convex optimization techniques \cite{bova04}.
Clearly, if $\sigma = 0$ then we
are back to the original Basis Pursuit principle (\ref{BP}).

For the problem (\ref{P2}) quite general stability results were obtained
by Candes, Romberg and Tao in \cite{carota06-1}, see also \cite{codade06}. 
Their key concept is the following definition. 

\begin{definition} The restricted isometry constant $\delta_M$ of a matrix
$A$ is the smallest number such that for all subsets $T$ with $|T| \leq M$
it holds
\begin{equation}\label{cond_RIP}
(1-\delta_M) \|x\|_2^2 \, \leq \, \|A_T x\|_2^2 \,\leq\, (1+\delta_M) \|x\|_2^2
\end{equation}
for all coefficients $x$ supported on $T$. Here $A_T$ denotes the restriction
of $A$ to the columns indexed by $T$.
\end{definition}

In \cite{carota06-1} the following theorem was proved. (Although it 
was originally stated
only for the real-valued case the theorem together with its proof 
also holds for the complex-valued case.)

\begin{Theorem}\label{thm_bp_stab} Assume that $A$ is some matrix for which
the restricted isometry constants satisfy
\[
\delta_{3M} + 3 \delta_{4M} < 2.
\]
Let $x \in \C^D$ and assume we have given noisy data $y=Ax+\eta$ with $\|\eta\|_2 \leq \sigma$. 
Denote by $x_M$ the truncated vector corresponding
to the $M$ largest absolute values of $x$. Then the solution $x^\#$
to the problem
\[
\min \|x\|_1 \mbox{ subject to } \|Ax-y\|_2 \leq \sigma
\]
satisfies
\begin{equation}\label{BP_error}
\|x^\# - x\|_2 \,\leq\, C_1 \sigma + C_2 \frac{\|x-x_M\|_1}{\sqrt{M}}.
\end{equation}
The constants $C_1$ and $C_2$ depend only on $\delta_{3M}$ and $\delta_{4M}$. 
\end{Theorem}

Thus, recovery by the BP variant (\ref{P2}) is stable provided the restricted
isometry constants are small. Note that the second 
term in (\ref{BP_error}) vanishes if $x$ is sparse, i.e.,
has not more than $M$ non-vanishing coefficients.

For our case this means that it is sufficient to provide conditions that
ensure $\delta_{4M} \leq \delta$ for some small $\delta$ 
with high probability. 
(Note that for $\delta = 1/5$ the constants in the previous
theorem are actually quite well-behaved, $C_1 \leq 12.04$ and $C_2 \leq 8.77$, see \cite{carota06-1}.) 

Cand{\`e}s and Tao \cite{cata06} provided such conditions 
for the discrete Fourier transform with a slightly different probability 
model than our discrete model (2).
More recently, Rudelson and Vershynin came up with a more elegant and shorter 
solution to 
this problem \cite{ru06-1}. 
It is possible to apply their technique also to our 
two probability models,
notably the continuous one. This gives the following result.

\begin{Theorem}\label{thm_BP} Let $D=|\Gamma|$ and a sparsity $M$ be given.
Let $\epsilon, \delta \in (0,1)$ and
assume
\begin{equation}\label{N_cond}
\frac{N}{\log(N)} \,\geq\, C \delta^{-2} M  \log^2(M) \log(D) \log(\epsilon^{-1}). 
\end{equation}
Let the $N$ sampling points $X=(x_1,\hdots,x_N)$ be chosen at random according
to the model (1) or (2). Then with probability 
at least $1-\epsilon$ the isometry constant of the matrix
$N^{-1/2} \F_X$ satisfies
\begin{equation}\label{delta}
 \delta_{M} \leq \delta.
\end{equation}
The constant $C$ is absolute.
\end{Theorem}

The combination of Theorems \ref{thm_bp_stab} and \ref{thm_BP} gives
the following.

\begin{corollary}\label{corBP} Let $\Gamma$ with $|\Gamma| = D$, 
$M$, $N$ and $\epsilon$ such that %be given as in
\begin{equation}\label{N_cond2}
\frac{N}{\log(N)} \,\geq\, C_0 M  \log^2(M) \log(D) \log(\epsilon^{-1}). 
\end{equation}
Choose $x_1,\hdots,x_N$ according to the probability
model (1) or (2). Then with probability at least $1-\epsilon$ the 
following holds for all coefficient vectors $c \in \C^\Gamma$. 
Assume $y = \F_X c+ \eta$ with $\|\eta\|_2 \leq \sigma$. Denote
by $c_M$ the truncated vector corresponding to the largest coefficients
of $c$. Then the solution
$c^\#$ to the minimization problem (\ref{P2})
satisfies
\begin{equation}\label{BP_approx_error}
\|c^\# - c\|_2 \leq C_1 \frac{\sigma}{\sqrt{N}} + C_2 \frac{\|c-c_M\|_1}{\sqrt{M}}.
\end{equation}
\end{corollary}

\begin{remark} \begin{itemize}
\item[(a)] Choosing $\sigma = 0$ yields uniform exact recovery. 
Under condition (\ref{N_cond}) 
BP is able to reconstruct exactly {\em all} $f \in \Pi_\Gamma(M)$ from
a single sampling set $X$.
\item[(b)] Note that condition (\ref{N_cond}) is satisfied if
$N \geq C \delta^{-2} M \log^4(D) \log(\epsilon^{-1})$. Furthermore, 
(\ref{N_cond2}) is probably not optimal. One may 
conjecture that $N = {\cal O}(M \log(D/\epsilon))$ or even 
$N={\cal O}(M\log(D/(M\epsilon)))$
samples are enough, see also \cite{ru06-1}. 
\item[(c)] With a discrete probability model 
(the variant of (2) outlined in Secion \ref{Sec_Prior}), 
Cand\`es and Tao originally obtained a version of Theorem \ref{thm_BP} 
(see \cite[Lemma 4.3]{cata06})
where for some parameter $\alpha$ and constant $\rho$ 
the statement $\delta_M \leq c_0$ holds
with probability at least $1 - C D^{-\rho/\alpha}$ under
the condition $N \geq \alpha^{-1} M \log(D)^6$.
Substituting $\epsilon = CD^{-\rho/\alpha}$ and solving for $\alpha$ 
yields the condition
\begin{equation}\label{CT_estim}
N \,\geq\, C' M \log(D)^5 \log(\epsilon^{-1}).
\end{equation}
It might be possible to adapt the original proof of Cand{\`e}s and
Tao also to the continuous probability model (1) although this
does not seem straightforward.
\end{itemize}
\end{remark}

\section{Orthogonal Matching Pursuit}

In this section we consider the stability of OMP. 
Since we measure only noisy samples we cannot expect
to have perfect recovery of a sparse signal, but at least
we would like to obtain the true support of the sparse
coefficient vector and only small deviations of their
entries. 
We first provide the analogue of Theorem \ref{thm:OMP_old} for the noisy
case. Unfortunately, we again have to restrict to the first iteration
because it is still not clear how to deal with the subtle stochastic
dependency issues arising in the analysis of the further iterations.

\begin{Theorem}\label{thm:OMP} Let $f \in \Pi_\Gamma(M)$ with
Fourier coefficients $c$. Let $N \in \N$ and $\tau,\epsilon \in (0,1)$ such
that
\begin{equation}\label{cond:samples_OMP}
N \geq C M \tau^{-2} \log(D / \epsilon).
\end{equation}
Further, let $\sigma > 0$ 
such that
\begin{equation}\label{cond:sigma}
\sigma \leq \frac{1-\tau}{4}\sqrt{\frac{N}{M}} \|c\|_2.
\end{equation}
Choose the random sampling set $X = (x_1,\hdots,x_N)$ according
to the probability model (1) or (2). Assume that we have given
noisy samples $y = (f(x_\ell))_{\ell = 1}^N + \eta = \F_X c + \eta$ 
with $\|\eta\|_2 \leq \sigma$.
Then with probability exceeding $1-\epsilon$ 
OMP selects an element of the true
support of $c$ in the first step. \\
If after $M$ steps OMP actually recovers the complete 
support of $c$ then with probability exceeding $1-\epsilon$ the reconstructed
coefficients $\tilde{c}$ satisfy
\begin{equation}\label{rec_error}
\|c -\tilde{c}\|_2 \leq \sqrt{\frac{2}{N}}\, \sigma.
\end{equation}
\end{Theorem}
From the proof of this Theorem one can deduce more precise information
about the constant in 
condition (\ref{cond:samples_OMP}). Indeed, $N$ has to satisfy
the two conditions
\[
N \geq 17.88\, M \tau^{-2} \log(8D/\epsilon)\quad \mbox{and} \quad
\left\lfloor \frac{N}{12e M} \right\rfloor \geq \ln(2(1-1/(4e))^{-1} M/\epsilon).
\]

Note that $\|c\|_2 \geq \sqrt{M} \min_{j \in T} |c_j|$. Hence, condition
(\ref{cond:sigma}) is satisfied if
\begin{equation}\label{sigma_min}
\sigma \leq \frac{1-\tau}{4} \sqrt{N} \min_{j \in T} |c_j|.
\end{equation}
One expects that this condition %(\ref{sigma_min}) 
(with possibly a different constant)
is sufficient that OMP selects an element of the true support $T$ in
every step. Hence, the noise level should not exceed the minimal absolute
non-zero coefficient in order to have recovery of the correct support.

We note that our numerical experiments in Section \ref{Sec:Numerical} 
indicate that under condition (\ref{cond:samples_OMP}) 
OMP actually selects elements of the true support $T$ also 
in the further iterations and then (\ref{rec_error}) holds. 
However, we have not yet been able to carry through 
the corresponding theoretical analysis.

\subsection{A uniform result}

The result in the previous section is non-uniform. Let us state also
a uniform recovery result for OMP extending Theorem \ref{thm:OMP_uniform}
to the noisy situation.

\begin{Theorem}\label{thm:uniformOMP} Let the random sampling 
set $X=(x_1,\hdots,x_N)$ be chosen
according to one of our probability models. Let $\tau,\epsilon \in (0,1)$ and $\sigma > 0$. Assume that
\begin{equation}\label{cond:uniform}
N \geq C \tau^2 (2M-1)^2 \ln(4D'/\epsilon),
\end{equation}
where $D' = \#\{j-k: j,k \in \Gamma, j\neq k\} \leq D^2$.
Then with probability $1-\epsilon$ the following holds for all
$f \in \Pi_\Gamma(M)$ whose Fourier coefficients satisfy
\begin{equation}\label{cond:noise2}
\min_{k \in \supp c} |c_k| > \frac{2\sigma}{(1-\tau)\sqrt{N}}.
\end{equation}
If OMP is applied on the noisy samples 
$y = \F_X c + \eta$ with $\|\eta\|_2 \leq \sigma$, and stopped once
the residual satisfies $\|r_s\| \leq \sigma$ then the true support of
$c$ is recovered and the reconstructed coefficient vector $\tilde{c}$
satisfies
\[
\|c - \tilde{c}\|_2 \leq \frac{1}{\sqrt{N(1-\tau/2)}}\, \sigma.
\]
\end{Theorem}
The above result has the drawback that 
the number of samples required by (\ref{cond:uniform})
scales quadratically in the sparsity $M$ rather than linearly 
as in (\ref{cond:samples_OMP}). 
As in the noiseless case however, 
one cannot expect to come around the quadratic scaling
if one requires uniformity, i.e., recovery by OMP of 
{\em all} $f \in \Pi_\Gamma(M)$
from a single sampling set $X$. Up to perhaps the $\log$-factor condition
(\ref{cond:uniform}) seems then to be optimal, see \cite{ra07}.

In contrast, BP gives a uniform guarantee if the number of samples
is only linear in the sparsity up to some $\log$-factors, see Theorem 
\ref{thm_BP}. 
Thus, under this
requirement, BP seems to be the method of choice. However, for certain
applications it might be enough to have a non-uniform guarantee and then
OMP is a good alternative considering that it is usually significantly faster
and much easier to implement, see also Section \ref{Sec:Numerical}.

\section{Proof of Theorem \ref{thm_BP}}

We mainly follow the ideas in \cite{ru06-1}.
Condition (\ref{cond_RIP}) for $N^{-1/2} \F_{X}$ is equivalent
to
\[
\sup_{|T| \leq M} \|I_T - N^{-1} \F_{TX}^* \F_{TX}\|_{2 \to 2} 
\,=\, \delta_M,
\]
and we have to prove that this inequality holds for $\delta_M \leq \delta$ 
with high probability. We denote by $z_\ell \in \C^\Gamma$ the vector
\begin{equation}\label{def_z_ell}
z_\ell \,=\, (e^{-ik\cdot x_\ell})_{k \in \Gamma}
\end{equation}
and by $z_\ell^T$ its truncation to the index set $T \subset \Gamma$. 
For vectors $y,z$ we define a rank one operator by $(y\otimes z)(x) = \langle x,y\rangle z$.
We note that
\[
(z_\ell^T \otimes z_\ell^T)(c) \,=\, \langle c, z_\ell^T \rangle z_\ell^T 
\,=\, \left( \sum_{j\in T} c_j e^{i(j-k)\cdot x_\ell} \right)_{k \in T}. 
\]
Observe that we can write $\F_{TX}^* \F_{TX} \,=\, \sum_{\ell=1}^N z_\ell^T \otimes z_\ell^T$.
Thus, we have to show that 
\begin{equation}\label{essential}
\sup_{|T| \leq M} \left\|I_T - \frac{1}{N} \sum_{\ell = 1}^N z_\ell^T \otimes z_\ell^T\right\|_{2\to 2} 
\leq \delta
\end{equation}
with probability at least $1-\epsilon$. To this end we consider the expectation of the above expression.
Further, we introduce an auxiliary matrix norm,
\[
\on A \on \,=\, \on A \on_M \,:=\, \sup_{|T| \leq M} \|A_{T,T}\|_{2 \to 2}
\]
where $A_{T,T}$ denotes the submatrix of a matrix $A$ 
consisting of the columns and 
rows indexed by $T$. The left hand side of 
(\ref{essential}) can be written as
\[
X_N \,:=\, \sup_{|T|\leq M} \left\| I_T - \frac{1}{N}\sum_{\ell=1}^N z_\ell^T \otimes z_\ell^T\right\|_{2\to 2} 
\,=\, \on I - \frac{1}{N} \sum_{\ell=1}^N z_\ell \otimes z_\ell\on 
\,=\, \on \sum_{\ell=1}^N N^{-1}(I-z_\ell \otimes z_\ell) \on .
\]
The random matrices $N^{-1}(I-z_\ell \otimes z_\ell)$, $\ell=1,\hdots,N$, 
are stochastically independent.
Moreover, it is easy to see that for both probability models
(1) and (2) $\E [z_\ell \otimes z_\ell] = I$ and 
$I - z_\ell\otimes z_\ell$ is symmetric. 
Then by standard symmetrization techniques,
see e.g.~\cite[Lemma 6.3]{leta91}, we have
\begin{align}
\E X_N &\,=\, \E \sup_{|T| \leq M} \left\|I_T - \frac{1}{N} \sum_{\ell = 1}^N z_\ell^T \otimes z_\ell^T\right\|_{2\to 2} 
\,=\, \E \left[\on I - \frac{1}{N}\sum_{\ell=1}^N z_\ell \otimes z_\ell \on\right]\notag\\
&\leq\, 2 \E \left[\on \frac{1}{N}\sum_{\ell=1}^N \epsilon_\ell z_\ell \otimes z_\ell \on \right]
\label{estimE1}
\,=\, 2 \E \sup_{|T| \leq M} \left\| \frac{1}{N} \sum_{\ell=1}^N \epsilon_\ell\, z_\ell^T \otimes z_\ell^T \right\|_{2\to 2},
\end{align}
where the $\epsilon_\ell$ are independent symmetric random variables taking values in $\{-1,+1\}$,
also jointly independent of the $x_\ell$.
Now the core of the proof is the following lemma
due to Rudelson and Vershynin \cite[Lemma 3.5]{ru06-1}. 
\begin{lemma}\label{ruve} Let $z_1,\hdots,z_N$, $N\leq D$, 
be (fixed) vectors in $\C^D$ with uniformly bounded entries, $\|z_\ell\|_\infty \leq 1$.
Then
\[
\E \sup_{|T| \leq M} \left\|\sum_{\ell=1}^N \epsilon_\ell\, z_\ell^T \otimes z_\ell^T \right \|_{2\to 2} 
\leq K(M,N,D) \sup_{|T| \leq M} \left\| \sum_{\ell=1}^N z_\ell^T \otimes z_\ell^T \right\|_{2 \to 2}^{1/2}
\]
where
\[
K(M,N,D) \,=\, C_0 \sqrt{M} \log(M) \sqrt{\log(D)}\sqrt{\log(N)}.
\] 
\end{lemma}
We remark that the elegant proof of this lemma uses entropy methods, 
in particular, Dudley's inequality \cite[Theorem 11.17]{leta91} for
the maximum of a Gaussian process.

Now, as in \cite{ru06-1}, we 
denote $E = \E[X_N]$. 
Using 
(\ref{estimE1}), taking the expectation only with respect to the
variables $\epsilon_\ell$, applying Lemma \ref{ruve} and H\"older's inequality 
we obtain
\begin{align}
E &\leq\, \frac{2 K(M,N,D)}{\sqrt{N}} 
\E \sup_{|T|\leq M} \left\|\frac{1}{N} \sum_{\ell=1}^N z_\ell^T \otimes z_\ell^T\right\|_{2\to2}^{1/2} \notag\\
&\leq\,  \frac{2 K(M,N,D)}{\sqrt{N}}\left(\E \sup_{|T|\leq M} \left\|I_T - \frac{1}{N} \sum_{\ell=1}^N z_\ell^T \otimes z_\ell^T\right\|_{2\to 2} + 1\right)^{1/2}
\,=\, \frac{2 K(M,N,D)}{\sqrt{N}} \sqrt{E+1}.\notag
\end{align}
It follows that $E\leq \theta$ provided 
\begin{equation}\label{condMND}
\frac{2 K(M,N,D)}{\sqrt{N}} \leq \frac{\theta}{\sqrt{1+\theta}}.
\end{equation}

To finish the proof we need to show that the random variable on the left hand 
side of (\ref{essential}) does not deviate much from its expectation. 
Inspired by \cite{caro07} we proceed differently as in \cite{ruve06} and use 
the following version of Talagrand's concentration inequality \cite{ta96-2}
proved by Klein and Rio in \cite{klri05}.

\begin{Theorem}\label{thm_Talagrand} Let $Y_1,\hdots,Y_N$ be a sequence of independent 
random variables with values in some Polish space $X$. Let $\cal{F}$ 
be a countable collection of real-valued
measurable and bounded functions $f$ on $X$ with $\|f\|_\infty \leq B$ for all 
$f \in \cal{F}$. Let $Z$ be the random variable
\[
Z = \sup_{f \in {\cal F}} \sum_{\ell=1}^N f(Y_\ell).
\]
Assume $\E f(Y_\ell) = 0$ for all $\ell=1,\hdots,N$ and all $f\in \F$.
Let $\sigma^2 := \sup_{f \in \F} \sum_{\ell=1}^N \E f(Y_\ell)^2$. Then for 
$t \geq 0$
\[
\P( Z \geq \E Z + t) \leq \exp\big(-\frac{t}{4B}\log\big(1+ 2\log\big(1+\frac{Bt}{2 B \E Z + \sigma^2}\big)\big)\big).
\]
\end{Theorem}

In order to apply the Theorem, we observe that
\begin{align}
X_N &= \sup_{|T|\leq M} \|I_T - \frac{1}{N}\sum_{\ell=1}^N z_\ell^T \otimes z_\ell^T\|
\notag\\
&=\, \sup_{|T| \leq M} \sup_{v \in \C^T, \|v\|_2 \leq 1}~\sup_{w \in \C^T, \|w\|_2 \leq 1}
\big| \frac{1}{N} \sum_{\ell=1}^N \langle (I_T - z_\ell^T \otimes z_\ell^T) v, w \rangle\big|
\notag\\
& = \sup_{(v,w) \in S_M} \big|\frac{1}{N} \sum_{\ell=1}^N \langle (I - z\otimes z) v, w\rangle\big|,\notag
\end{align}
where 
\[
S_M = \{(v,w) \in \C^D: \|v\|_2,\|w\|_2 \leq 1; \supp v = \supp w = T
\mbox { for some } T \mbox{ with } |T| \leq M\}.
\]
Defining
\[
f_{v,w}(z) \,=\, \frac{1}{N} \langle (I-z \otimes z) v, w\rangle.
\]
we obtain
\[
X_N = \sup_{(v,w) \in S_M} |\sum_{\ell=1}^N f_{v,w}(z_\ell)|.
\]
Clearly, $\E f_{v,w}(z_\ell) = N^{-1} \langle \E(I-z_\ell \otimes z_\ell) v, w\rangle = 0$.
Furthermore, for $(v,w) \in S_M$ and $z = (e^{ik\cdot x})_{k \in \Gamma}$ we have
\begin{align}
|f_{v,w}(z)| &= \frac{1}{N} \big|\sum_{j,k \in T, j\neq k} v_j e^{i(j-k)x} \overline{w_k}\big|
 \leq N^{-1} \sum_{j,k \in T, j\neq k} |v_j| |w_k| \notag\\
& \leq  N^{-1} \sum_{j \in T} \sum_{k\in T} |v_k| |w_{\sigma_j(k)}|
= N^{-1} \sum_{j \in T} \langle |v|, |w^{(\sigma_j)}|\rangle 
\leq N^{-1} M \|v\|_2 \|w\|_2
\leq \frac{M}{N},\label{vw_est}
\end{align}
where $\{(k,\sigma_j(k)), j,k \in T\}$ is a reparametrization 
of $T \times T$ such that $\sigma_j(T) = T$, i.e., $\sigma_j$ is a 
suitable permutation; and $w^{(\sigma_j)}$ denotes the corresponding 
vector of reordered entries of $w$.
Above we used the Cauchy Schwarz inequality in the
fifth step. We deduced $\|f_{v,w}\|_\infty \leq M/N$ for all $(v,w) \in S_M$.

Next, for $(v,w) \in S_M$, we compute
\begin{align}
\E |f_{v,w}(z_\ell)|^2 \,&=\, N^{-2}\E \big|\sum_{j,k \in T, j\neq k} v_j e^{i(j-k)\cdot x_\ell} \overline{w_k}\big|^2\notag\\
\,&=\, N^{-2} \sum_{j,k \in T, j\neq k}~ \sum_{j',k' \in T, j'\neq k'} v_j \overline{v_{j'}}\, \overline{w_k} w_{k'} \E[e^{i(j-k-j'+k')\cdot x_\ell}]. \notag  
\end{align}
Since $x$ is uniformly distributed on $[0,2\pi]^d$ or on $\frac{2\pi}{m} \Z^d_m$
we have
$\E[e^{i(j-k-j'+k')\cdot x_\ell}] = \delta_{j',j-k+k'}$ and, hence,
\begin{align}
\E |f_{v,w}(z_\ell)|^2 \,&=\, N^{-2}\sum_{j,k \in T, j\neq k}~ \sum_{k' \in T} v_j \overline{v_{j-k+k'}}
w_{k'} \overline{w_{k}} 
\leq N^{-2}\|v\|_2^2 \sum_{k,k' \in T} w_{k'} \overline{w_k}\notag\\
&\leq N^{-2} \|v\|_2^2\, |T|\, \|w\|_2^2 
\leq M/N^2.\notag
\end{align}
In the second step we applied the Cauchy Schwarz inequality and in the
third step a similar estimate as in (\ref{vw_est}).
Hence,
\[
\sigma^2 = \sup_{(v,w) \in S_M} \sum_{\ell=1}^N \E |f_{(v,w)}(z_\ell)|^2 \leq M/N.
\]
Theorem 5.2 applies to real-valued functions $f$. Hence, we split into
real and imaginary parts $f_{v,w}^r = \operatorname{Re}(f_{v,w})$, 
$f_{v,w}^i =  \operatorname{Im}(f_{v,w})$.
Then the estimates above apply also to these functions, i.e.,
$\|f_{v,w}^r\|_\infty, \|f_{v,w}^i\|_\infty \leq \frac{M}{N}$ and
$\sigma_r^2, \sigma_i^2 \leq \frac{M}{N}$. 

Denote $Z^r = \sup_{(v,w) \in S_M} \sum_{\ell=1}^N f_{v,w}^r(z_\ell)$ and similarly
define $Z^i$. Since $f_{v,-w} = -f_{v,w}$ we have 
$Z^r = \sup_{(v,w) \in S_M} |\sum_{\ell=1}^N f_{v,w}(z_\ell)|$. 
By the union bound
\[
\P (X_N \geq \delta) = \P(|Z^r|^2 + |Z^i|^2 \geq \delta^2)
\leq \P\big(Z^r \geq \frac{\delta}{\sqrt{2}}\big) + \P\big(Z^i \geq \frac{\delta}{\sqrt{2}}\big).
\] 
Now assume $\E X_N \leq \delta/2$, which by (\ref{condMND}) will be satisfied 
provided $2K(M,N,D)/\sqrt{N} \leq \frac{\delta/2}{\sqrt{1 + \delta/2}}$, in particular, if 
\[
2K(M,N,D)/\sqrt{N} \leq \frac{\delta}{2\sqrt{3/2}} = \frac{\delta}{\sqrt{6}}.
\]
Setting $t=\frac{\delta}{\sqrt{2}} - \frac{\delta}{2} = \frac{\sqrt{2}-1}{2}\delta =: c \delta$ in Theorem \ref{thm_Talagrand} we obtain
\[
\P(X_N \geq \delta) \leq 2 \exp\left(-\frac{t}{4M/N} \log(1 + 2\log(1 + \frac{t}{\delta + 1}))\right) = 2 e^{-c_0(\delta) N/M},
\]
where $c_0(\delta) = c\delta \log(1+2\log(1+ \frac{c\delta}{\delta + 1}))$.
In other words, $X_N \leq \delta$ with probability at least $1-\epsilon$ provided
$N \geq c_0(\delta)^{-1} M \log(2/\epsilon)$ and $2 K(M,N,D)/\sqrt{N} \leq \delta/\sqrt{6}$.
Note that $C_2 \delta^{-2} \geq c_0(\delta)^{-1}$ for all $\delta \in (0,1)$ where
$C_2 = c_0(1)^{-1} = \frac{2}{\sqrt{2}-1} \log(1+2\log(1+\frac{\sqrt{2}-1}{4}))^{-1} \approx 26.84$.
With the definition of $K(M,N,D)$ %and since $c_0(\delta) \leq c \delta^2$ 
we deduce that $\delta_M \leq \delta$
with probability at least $1-\epsilon$ provided
\[
\frac{N}{\log(N)} \geq C_1 \delta^{-2} M \log^2(M) \log(D)
\quad \mbox{ and } \quad N \geq C_2 \delta^{-2} M \log(2\epsilon^{-1}).
\]
Both conditions are satisfied once
\[
\frac{N}{\log(N)} \geq C \delta^{-2} M \log^2(M) \log(D) \log(\epsilon^{-1})
\]
for some suitable constant $C$. This finishes the proof of Theorem \ref{thm_BP}.

\section{Proofs for Orthogonal Matching Pursuit}
\label{sec_Proof_OMP}

\subsection{Proof of Theorem \ref{thm:OMP}}

The proof is an extension of the one in \cite{kura06}.
We will use the following result from \cite{grpora07} on
the eigenvalues of a submatrix $\F_{TX}$, which is based on the analysis
in \cite[Lemma 3.3 and Section 3.3]{ra05-7}.

\begin{Theorem}\label{thm:eigvals} 
Let $T$ of
size $|T| = M$ and let $x_1,\hdots,x_N$ be
i.i.d. random variables that are uniformly distributed over
$[0,2\pi]^d$ or over the grid $\frac{2\pi}{m} \Z_m^d$.
Choose $\epsilon, \delta \in (0,1)$ and assume
\begin{equation}\label{cond:eigvals}
\left \lfloor \frac{\delta^2 N}{3e M} \right\rfloor \,\geq\,
\ln (c(\delta) M/\epsilon),
\end{equation}
where $c(\delta) = (1-\delta^2/e)^{-1} \leq (1-e^{-1})^{-1} \approx 1.582$.
Then with probability at least $1-\epsilon$ the minimal and
maximal eigenvalue
of $\F_{TX}^* \F_{TX}$ satisfy 
\begin{equation}\label{eigs}
1-\delta \leq \lambda_{\min}(N^{-1} \F_{TX}^* \F_{TX}),
\quad \mbox{and} \quad \lambda_{\max}(N^{-1} \F_{TX}^* \F_{TX})
\leq 1+\delta.
\end{equation}
\end{Theorem}

Further, we need the following concentration inequality proved in 
\cite{kura06}.

\begin{lemma}\label{lem:conc}
Assume that $c$ is a vector supported on $T$. 
Further, assume that the sampling set $X$ is chosen according
to one of our two probability models. 
Then for $j \notin T$ and $t > 0$ it holds
\[
\P\left(|N^{-1} \langle \F_{TX} c, \phi_j \rangle | \geq t\right)
\leq 4\exp\left(- N \frac{t^2}{4\|c\|_2^2 + \frac{4}{3\sqrt{2}}\|c\|_1 t}\right).
\] 
\end{lemma}

Now we can turn to the proof of Theorem \ref{thm:OMP}. (Orthogonal)  Matching
Pursuit selects an element of the support $\supp c =: T$ in the first iteration if
\begin{equation}\label{rec:cond}
\max_{j \notin T} 
|N^{-1} \langle \phi_j, \F_{TX} c + \eta\rangle| < 
\max_{k \in T} |N^{-1} \langle \phi_k, \F_{TX} c +\eta \rangle |.
\end{equation}
By the triangle inequality and Cauchy-Schwarz 
(note also that $\|\phi_k\|_2 = \sqrt{N}$) this will be satisfied if
\[
\max_{j \notin T} |N^{-1} \langle \phi_j, \F_{TX} c \rangle|
\leq \| N^{-1} \F_{TX}^* \F_{TX} c\|_\infty - \frac{2}{\sqrt{N}} \|\eta\|_2.
\]
Assume for the moment that $\lambda_{\min}(N^{-1} \F_{TX}^* \F_{TX}) 
\geq 1-\delta$ for some $\delta \in (0,1)$. (The probability that this happens
can be estimated by Theorem \ref{thm:eigvals}.) This yields
\[
\| N^{-1} \F_{TX}^* \F_{TX} c\|_\infty \geq M^{-1/2} 
\|N^{-1} \F_{TX}^*\F_{TX} c \|_2 \geq M^{-1/2} (1-\delta) \|c\|_2.
\]
Thus, (\ref{rec:cond}) is satisfied if
\begin{equation}\label{rec:cond2}
\max_{j \notin T} |N^{-1} \langle \phi_j, \F_{TX} c \rangle|
\leq \frac{1-\delta}{\sqrt{M}} \|c\|_2 - \frac{2}{\sqrt{N}} \|\eta\|_2.
\end{equation}   
Assuming further that
\begin{equation}\label{cond:epsilon}
\|\eta\|_2 \leq \frac{(1-\delta)(1-\tau)}{2}\sqrt{\frac{N}{M}} \|c\|_2
\end{equation}
condition (\ref{rec:cond2}) becomes true if
\[
\max_{j \notin T} |N^{-1} \langle \phi_j, \F_{TX} c \rangle|
\leq \frac{1-\delta}{\sqrt{M}} \tau \|c\|_2.
\]
By the concentration inequality in Lemma \ref{lem:conc} the probability
that the above inequality does not hold can be estimated by
\begin{align}
& \P\left(\max_{j \notin T} |N^{-1} \langle \phi_j, \F_{TX} c \rangle|
\leq \frac{1-\delta}{\sqrt{M}} \tau \|c\|_2 \right)%\notag\\
\,\leq\, \sum_{ j \notin T} \P\left( |N^{-1} \langle \phi_j, \F_{TX} c \rangle|
\leq \frac{1-\delta}{\sqrt{M}} \tau \|c\|_2 \right)\notag\\
&\leq 4 D \exp\left(-\frac{N}{M} \frac{(1-\delta)^2 \tau^2 \|c\|_2^2}
{4 \|c\|_2^2 + \frac{4}{3\sqrt{2}} \|c\|_1 M^{-1/2}(1-\delta)\tau\|c\|_2}  \right)\notag\\
&\leq 4D \exp\left(-\frac{N}{M} \frac{(1-\delta)^2 \tau^2}{4 + \frac{4}{3 \sqrt{2}}(1-\delta)}\right).   
\end{align}
In the last line we used the Cauchy-Schwarz inequality, 
$\|c\|_1 \leq \sqrt{M} \|c\|_2$. Now we choose $\delta = 1/2$. Then condition
(\ref{cond:epsilon}) becomes (\ref{cond:sigma}) and
\[
 \P\left(\max_{j \notin T} |N^{-1} \langle \phi_j, \F_{TX} c \rangle|
\leq \frac{1-\delta}{\sqrt{M}} \tau \|c\|_2 \right)
\leq 4D \exp\left(-\frac{N}{M} \frac{\tau^2}{16 + \frac{8}{3\sqrt{2}}}\right).
\]
The latter term is less than $\epsilon/2$ if
\[
 N \geq C M \tau^{-2} \log(8D/\epsilon)
\]
with $C = 16 + \frac{8}{3\sqrt{2}} \approx 17.88$. Furthermore, by Theorem \ref{thm:eigvals}
our initial assumption that $\lambda_{\min}(N^{-1} \F_{TX}* \F_{TX}) \geq 1-\delta = 1/2$ fails with probability at most $\epsilon/2$ if
\[
\left \lfloor \frac{N}{12 e M} \right \rfloor \geq \ln(2(1-e^{-1}/4)^{-1} M /\epsilon). 
\]
Altogether, the probability that OMP does not select an element of $T$ in the 
first step is less than $\epsilon$ if
\[
N \geq C M \tau^{-2} \log(D/\epsilon)
\]
for some suitable constant $C$. 

Now consider the final statement of the Theorem, i.e.,
assume that OMP has reconstructed the true support $T$ after $M$ steps.
Then the reconstructed coefficients are given by
$\tilde{c} = \F_{TX}^\dagger (\F_{TX} c + \eta)$ where $\F_{TX}^\dagger$ denotes
the pseudo-inverse of $\F_{TX}$. Observe that $\F_{TX}^\dagger \F_{TX} c = c$. 
Hence
\[
\|\tilde{c}-c\|_2 \,=\, \|\F_{TX}^\dagger \eta\|_2
\leq \on \F_{TX}^\dagger \on \|\eta\|_2 
= \sqrt{\lambda_{\min}(\F_{TX}^* \F_{TX})^{-1}} \|\eta\|_2 
 \leq  \frac{1}{\sqrt{N(1-\delta)}}\, \sigma 
= \sqrt{\frac{2}{N}}\, \sigma, 
\]
where we used Theorem \ref{thm:eigvals} once more.

\subsection{Proof of Theorem \ref{thm:uniformOMP}}

The proof of the uniform recovery result is based on the 
coherence parameter, which measures the maximum correlation
between distinct columns of a matrix $A= (\psi_1|\hdots|\psi_D)$, i.e.,
\[
\mu \,:=\,\max_{j \neq k} |\langle \psi_j,\psi_k \rangle|.
\]
Based on $\mu$ the following theorem
due to Donoho, Elad and Temlyakov \cite[Theorem 4.1]{doelte06} 
analyzes the performance
of OMP in the presence of noise.

\begin{Theorem}\label{thm:doelte} Assume that $A$ has 
coherence $\mu$. Suppose
that $y = A c + \eta$ with only $M$ coefficients of $c$ being nonzero
and $\|\eta\|_2 \leq \sigma$. Suppose $(2M-1)\mu < 1$ and
\[
\sigma < \frac{1-(2M-1)\mu}{2} \min_{k \in \supp c} |c_k|.
\]
If we run OMP until the residual satisfies $\|r_s\| \leq \sigma$ then
the true support of $c$ has been recovered, and consequently OMP has done
$M$ iterations. Furthermore, the error between the reconstructed
coefficients $\tilde{c}$ and the original coefficients satisfies
\[
\|c-\tilde{c}\|_2 \leq \frac{1}{\sqrt{1-(M-1)\mu}}\, \sigma.
\]
\end{Theorem}

In \cite{kura06} the following estimate of the coherence of $\F_X$ was
proven.

\begin{lemma} Let the random sampling set $X= (x_1,\hdots,x_N)$ be chosen
according to one of our probability models and let $\mu$ be the
coherence of the random matrix $N^{-1/2}\F_{X}$. Then
\[
\P(\mu > t) \leq 4 D' \exp\left(-N \frac{t^2}{4+\frac{4}{3\sqrt{2}}t}\right),
\]
where $D' = \#\{j-k: j,k \in \Gamma, j\neq k\} \leq D^2$.
\end{lemma}
\begin{remark} In case of the continuous probability model
the previous estimate can be slightly improved to \cite{kura06} 
\[
\P(\mu > t) \leq (1-\kappa)^{-1} D' e^{-N\kappa t^2},\qquad \kappa \in (0,1).
\] 
\end{remark}
Now the proof of Theorem \ref{thm:uniformOMP} is a mere application
of the above statements. Note 
that $y = \F_X c + \eta = N^{-1/2} \F_X (\sqrt{N} c) + \eta$. 
Thus, setting $t = \frac{\tau}{2M-1}$ in the lemma, solving for $N$ and
using Theorem \ref{thm:doelte} with $c' = \sqrt{N} c$ 
shows the assertion.

\section{Numerical Experiments}
\label{Sec:Numerical}

\begin{figure}%[h]
  \centering
  \subfigure[Recovery success rate of the true support.]   
  {\includegraphics[width=0.45\textwidth]{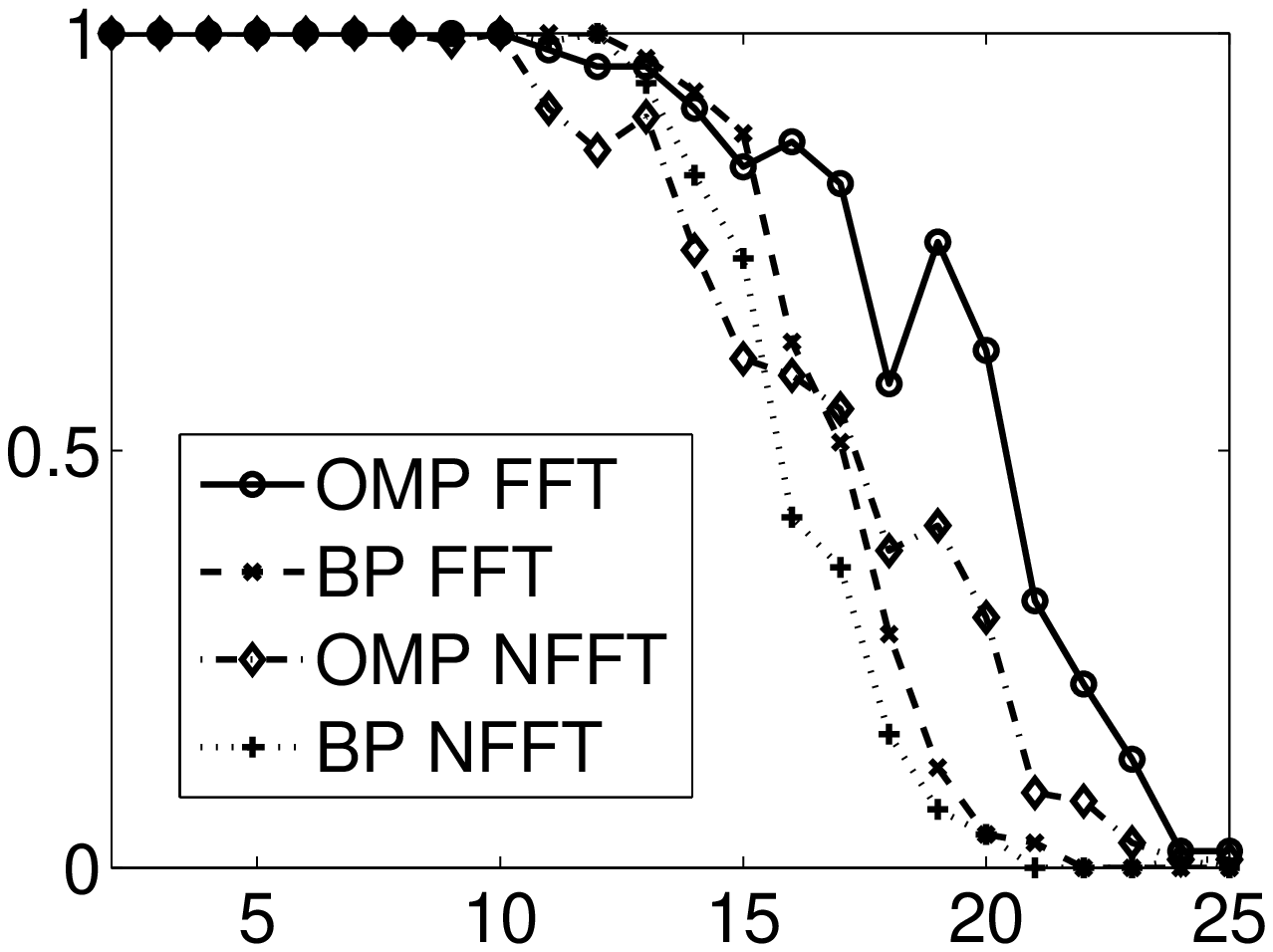}}\hfill
  \subfigure[Avarage error.]
  {\includegraphics[width=0.45\textwidth]{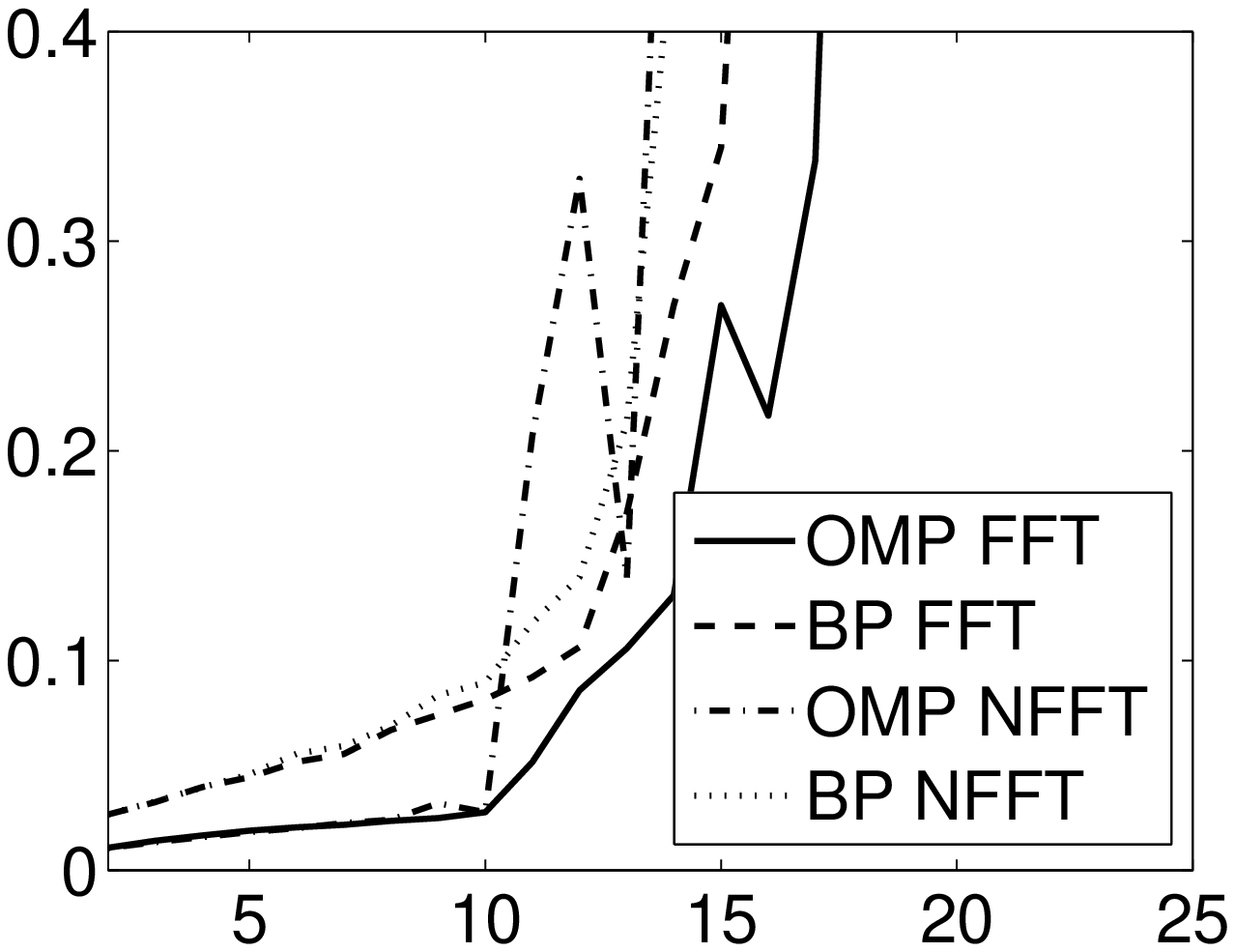}}\\  
  \caption{Recovery of sparse trigonometric polynomials in dimension $D=256$ 
from $N=50$ noisy samples, $\|\eta\|_2 = \sigma = 0.4$. 
The sparsity $M$ is varied.\label{fig:sigma04}}
\end{figure}

To illustrate the theoretical results we also conducted 
numerical experiments. We choose a number
of samples $N$, the noise level $\sigma$, the sparsity $M$ and an 
(even) dimension $D$ and set $\Gamma = \{-D/2+1,\hdots,D/2\}$. 
Then we repeat the following reconstruction experiment $100$ times.
We choose a subset $T$ uniformly at random among all
subsets of size $M$. Then we randomly 
select the real part and imaginary part of the 
coefficients $c_k$ on $T$ from a standard normal distribution.  
The sampling points $x_1,\hdots,x_N$ are randomly drawn either from the
uniform distribution on $[0,2\pi]$ (probability model (1), labelled 
NFFT in the plots) 
or uniformly among all subsets of 
$\{0,\frac{2\pi}{D},\hdots,\frac{2\pi(D-1)}{D}\}$ of size $N$ 
(a slight variation of the probability model (2) 
preventing that some of the sampling points coincide, labelled FFT). 
The perturbed sampling points
are given by $y_\ell = \sum_{k \in T} c_k e^{ik\cdot x_\ell} + \eta_\ell$, 
$\ell = 1,\hdots,N$, where the noise vector $\eta$ is chosen
uniformly at random on the sphere with radius $\sigma$ in $\C^N$, 
i.e.~$\|\eta\|_2 = \sigma$.

Then we solve the $\ell_1$-minimization problem (\ref{P2}) 
(with the chosen $\sigma$) 
and run OMP (with precisely $M$ iterations), respectively, and
compute the error between the reconstructed vector $\tilde{c}$ and 
the original vector $c$ for both methods. 
Also we test whether the correct support has been recovered.

\begin{figure}%[h]
  \centering
\subfigure[Average error, $\sigma = 5$, $D=256$, $N=50$]   
  {\includegraphics[width=0.45\textwidth]{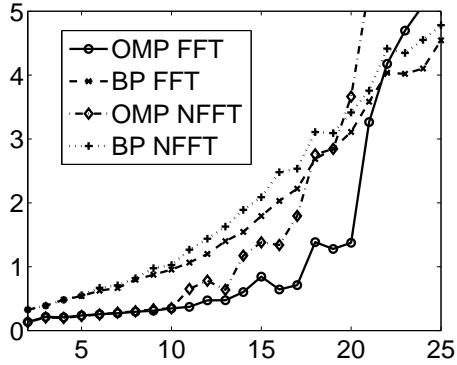}}\hfill
  \subfigure[Average error, $\sigma = 2$, $D=512$, $N=150$]
  {\includegraphics[width=0.45\textwidth]{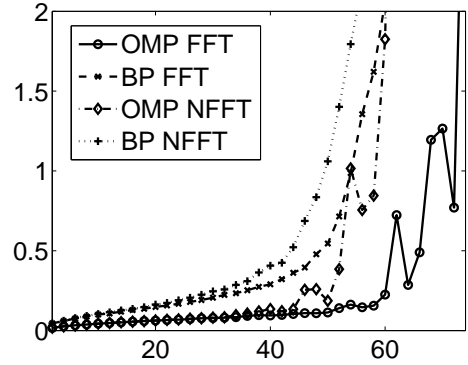}}\\
  \caption{Recovery of sparse trigonometric polynomials for different
sets of parameters. The sparsity is varied.\label{fig:sigma2}}
\end{figure}

\begin{figure}%[h]
  \centering
\subfigure[Recovery success, $D=256$, $N=50$, $M=10$]   
  {\includegraphics[width=0.45\textwidth]{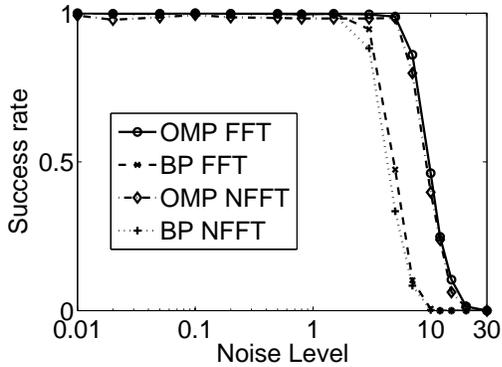}}\hfill
  \subfigure[Average error, $D=256$, $N=50$, $M=10$]
  {\includegraphics[width=0.45\textwidth]{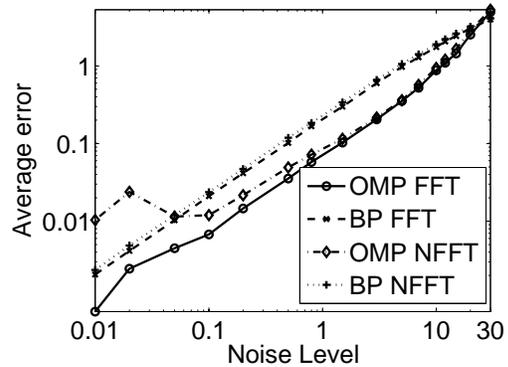}}\\
  \caption{Recovery of sparse trigonometric polynomials for different
sets of parameters. The noise level $\sigma$ is varied.
(For comparison, the average $\ell_2$-norm of the vector of samples 
of the unperturbed polynomial is approximately 39.4). 
\label{fig:sigma5}}
\end{figure}

Figures \ref{fig:sigma04} and \ref{fig:sigma2} show the results
for varying sparsity, while in Figure \ref{fig:sigma5} the noise level
$\sigma$ is varied.  
These plots indicate that the BP variant and OMP are both stable under noise 
as predicted by the theoretical results. Figure 4 suggests that the correct
support set can be recovered even when the noise level reaches the order
of the $\ell_2$-energy of the samples of the signal.
Moreover, OMP usually 
performs slightly better than BP.
In fact, OMP yields a smaller avarage reconstruction error and also
reconstructs more often the correct support - despite
that fact that theoretically BP gives a uniform recovery guarantee 
while OMP does not. This might be due to the fact that OMP forces the 
reconstruction to be $M$-sparse while BP may result in larger support sets.
Furthermore, OMP
is much faster than BP (by a factor between $10$ and $200$ in the 
examples). For a more detailed comparison of the computation times
we refer to \cite{kura06}. 

The Matlab toolbox CVX \cite{CVX} was used for solving
(\ref{P2}).
The examples (including the OMP algorithm) are part of 
the Matlab toolbox \cite{kura_OMP}, which is available online.

\section{Discussion}

We presented theoretical and numerical results concerning the stability 
of recovery of sparse trigonometric polynomials with (a variant
of) Basis Pursuit and Orthogonal Matching Pursuit.
The (non-uniform) recovery Theorem \ref{thm:OMP} for OMP, however,
is only partial so far. It remains open to analyze theoretically 
the further iterations
after the first step. 

BP has the advantage of giving a uniform guarantee of recovery
success, i.e., a single sampling set $X$ may be sufficient to recover
{\em all} sparse trigonometric polynomials, while it seems that OMP
is only able to provide non-uniform recovery results at reasonably small
ratio of the number of samples to the sparsity \cite{ra07}. (But note 
the results for variants of OMP in \cite{neve07,neve07-1,netr08}.)
In practice, however, a non-uniform guarantee might be sufficient
and indeed our numerical experiments show that OMP even slightly outperforms
BP on generic (=random) signals. 

Corollary \ref{corBP} concerning BP covers also the case that the coefficient
vector $c$ is not sparse in a strict sense. 
In this case it estimates the approximation 
error of the reconstruction by the approximation error with $M$-terms.
In principle, we might also apply Theorems \ref{thm:OMP} and 
\ref{thm:uniformOMP} for OMP 
to the non-sparse case by letting $\eta = \F_X c_{\Gamma \setminus T}$, i.e.,
by treating the contribution of the (small) coefficients outside $T$ as noise.
However, for most situations conditions (\ref{sigma_min})
and (\ref{cond:noise2}) on the magnitude of
the coefficients become then unrealistic. Roughly speaking they would imply 
that the smallest coefficient of $c$ in $T$ is significantly larger
than the $\ell_2$-norm of the coefficients outside $T$.
So a thorough treatment of the non-sparse case for OMP is still open. 

OMP is usually faster (and easier to implement) 
than BP in practice, and the numerical results even
indicate that OMP is slightly more stable. So in most practical situations 
one would probably prefer to use OMP despite its lack of giving
a uniform recovery guarantee when the number of samples is only linear
in the sparsity.

\section*{Acknowledgements}

The author would like to thank Stefan Kunis, Ingrid Daubechies,
Joel Tropp, R{\'e}mi Gribonval, Pierre Vandergheynst, Karin Schnass 
and Roman Vershynin for enlightening discussions. 
Parts of the manuscript were written
while visiting the group of Pierre Vandergheynst at EPFL in Lausanne. He thanks
its members for their warm hospitality. His stay was funded by the European
Union's Human Potential Program under contract HPRN-CT-2002-00285 (HASSIP).
The author is currently 
supported by an Individual Marie Curie fellowship from the 
European Union under contract MEIF-CT-2006-022811.

\bibliographystyle{abbrv}
\bibliography{RauhutBib}

\end{document}